\documentclass[amsmath,secnumarabic,floatfix,amssymb,nofootinbib,nobibnotes,letterpaper,11pt,tightenlines]{revtex4}

\usepackage{times}

\usepackage{geometry}
\usepackage{amssymb}
\usepackage{latexsym, amsmath, amscd,amsthm}
\usepackage{graphicx}
\usepackage[percent]{overpic}
\usepackage{units}
\usepackage{hyperref}
\PassOptionsToPackage{caption=false}{subfig}
\usepackage[lofdepth]{subfig}
\usepackage{algorithm}
\usepackage{algorithmicx}
\usepackage{algpseudocode}
\usepackage{mathtools}

\newtheorem{theorem}{Theorem}

\newtheorem{proposition}[theorem]{Proposition}

\theoremstyle{definition}

\usepackage{units}
\newcommand{\R}{\mathbb{R}}

\newcommand{\Z}{\mathbb{Z}}
\newcommand{\C}{\mathbb{C}}

\newcommand{\RP}{\mathbb{RP}}

\makeatletter
\renewcommand*\env@matrix[1][*\c@MaxMatrixCols c]{%
  \hskip -\arraycolsep
  \let\@ifnextchar\new@ifnextchar
  \array{#1}}
\makeatother

\begin{document}

\title[]{Spherical Geometry and the Least Symmetric Triangle}
\author{Laney Bowden}
\affiliation{Colorado State University, Fort Collins, CO}
\author{Andrea Haynes}
\affiliation{Colorado State University, Fort Collins, CO}
\author{Clayton Shonkwiler}
\affiliation{Colorado State University, Fort Collins, CO}
\author{Aaron Shukert}
\affiliation{Colorado State University, Fort Collins, CO}

\begin{abstract}
We study the problem of determining the least symmetric triangle, which arises both from pure geometry and from the study of molecular chirality in chemistry. Using the correspondence between planar $n$-gons and points in the Grassmannian of 2-planes in real $n$-space introduced by Hausmann and Knutson, this corresponds to finding the point in the fundamental domain of the hyperoctahedral group action on the Grassmannian which is furthest from the boundary, which we compute exactly. We also determine the least symmetric obtuse and acute triangles. These calculations provide prototypes for computations on polygon and shape spaces.
\end{abstract}

\maketitle

The equilateral triangle is surely the most symmetric triangle by any reasonable standard, but what is the \emph{least} symmetric triangle? More generally, what is the least symmetric $n$-gon?

As stated, this question is ill-posed, but we can make sense of it using a map defined by Hausmann and Knutson~\cite{Hausmann:1997p8571} from the Grassmannian $G_2(\R^n)$ of 2-dimensional linear subspaces of $\R^n$ to the collection of ordered planar $n$-gons up to similarity.\footnote{Here \emph{ordered} means that the order of edges in the $n$-gon matters: e.g., a triangle with edges ordered from shortest to longest is distinct from the same triangle with edges ordered from longest to shortest.} The symmetric group $S_n$ acts on ordered $n$-gons by permuting edges, and this action lifts to an action of the hyperoctahedral group $B_n \simeq (\Z/2\Z)^n \rtimes S_n$ on $G_2(\R^n)$ (see~\cite{random-triangles}). The fundamental domain of this action is a region in $G_2(\R^n)$ whose boundary corresponds to those $n$-gons which are unchanged by some group element: in other words, the boundary corresponds to $n$-gons with a symmetry.

Therefore, we can re-state the question as: which point(s) are furthest from the boundary of the fundamental domain of the hyperoctahedral group action on $G_2(\R^n)$? The answer should be important in understanding the behavior of permutation-invariant functions on random walks, and should also be interesting from the perspective of algebraic geometry. We like this question, but don't know the general answer. So in this paper, we pose the question precisely and work out the answer for triangles in as much detail as possible in order to jumpstart the larger problem.

In the case of triangles, the Grassmannian $G_2(\R^3)$ is double-covered by the unit sphere $S^2$, which is where we will actually do our calculations. In $S^2$, the fundamental domain is simply a spherical triangle that, as we will see, can be interpreted as the space of \emph{unordered} triangles. Its boundary is precisely the set of isosceles and degenerate triangles and its interior is the set of scalene triangles; note that both isosceles and degenerate triangles have a mirror symmetry. Therefore, the point maximizing the distance to the boundary -- which will be our least symmetric triangle -- could just as well be thought of as the \emph{most scalene triangle}. Our triangle will be different from any of those proposed by Robin~\cite{Robin:2009vx}.

The problem of finding the most scalene triangle also arises in chemistry, where the chirality of molecules plays a key role~\cite{Prelog:1976dk}. A shape is \emph{chiral} if it is not congruent to its mirror image; otherwise it is \emph{achiral}.

Although three-dimensional notions of chirality are the most physically relevant, substantial effort has been expended on the study of chirality in two dimensions, which ``serves as the first step for a deeper understanding of three-dimensional chirality in chemistry''~\cite{Buda:1991fp}. In the case of triangles, isosceles triangles are achiral since reflection across the median intersecting the odd edge produces a congruent triangle, while scalene triangles are chiral.\footnote{Strictly speaking, scalene triangles are only chiral when viewed as living in a two-dimensional universe. In three dimensions a scalene triangle and its mirror image are related by a rotation.} Consequently, there is an extensive chemistry literature on scalene triangles and in particular the search for the most scalene triangle~\cite{Buda:1991fp,Buda:1992kc,Chauvin:1996eh,Rassat:2003jn}.


In the chemistry literature, the most scalene (or most chiral) triangles are found as maxima of certain energy functions defined on various models of triangle space. These models are mostly \emph{ad hoc}, whereas we view the Grassmannian parametrization of polygon space as a principled choice of model: after all, for any $n$ it yields a model of $n$-gon space which is a Riemannian symmetric space, meaning in particular that it has a transitive group of isometries and, since it is compact, a canonical (up to scale) invariant Riemannian metric. 

Moreover, this model is beginning to gain traction in polymer physics due to its computational tractability~\cite{Suzuki:2014fo,Uehara:2014kf,Deguchi:2017bz}, and a version of it for continuous curves has been independently developed for use in shape recognition and classification problems~\cite{Younes:2008gy}. We see the current paper as a prototype of a geometric approach to finding optimal polygons or shapes.

\section{The construction}\label{sec:construction}

We introduced the \emph{symmetric measure} on the space of planar $n$-gons in~\cite{Cantarella:2013bla} by pushing forward the uniform measure on $G_2(\R^n)$ using Hausmann and Knutson's~\cite{Hausmann:1997p8571} map from the Grassmannian $G_2(\R^n)$ whose points correspond to the 2-dimensional linear subspaces of $\R^n$. In other words, a plane in $\R^n$ corresponds to a similarity class of $n$-gons in the plane, and the measure of a collection of polygons is the measure of the corresponding subset of the Grassmannian. In fact, we get much more than just a probability measure on polygon space: defining the Hausmann--Knutson map to be a Riemannian submersion yields a Riemannian metric on polygon space.

Something special happens when $n=3$: any 2-dimensional subspace of $\R^3$ has a unique line through the origin as its orthogonal complement, so we can equivalently think about the space of lines through the origin in $\R^3$, otherwise known as the projective plane $\RP^2$. In turn, the projective plane is double-covered by the unit sphere $S^2 \subseteq \R^3$ since any unit vector determines a line through the origin, and the only way that two distinct unit vectors can lie in the same line is if they are antipodal. Consequently, we can identify the space of planar triangles with the unit sphere in three-dimensional space.

This construction and identification (and the generalization to Grassmannians) is explained in~\cite{random-triangles}, but we briefly summarize it here. 

The space of planar triangles is non-compact since the size of any triangle can be scaled by an arbitrary positive real number. However, since we will primarily be interested in triangles up to similarity, we may as well compactify by fixing a scale for our triangles. As will shortly be apparent, it turns out to be convenient to normalize all triangles to have perimeter 2.

Since triangles are determined up to congruence by their three side lengths, we can uniquely specify a triangle up to similarity by choosing the side lengths $a$, $b$, $c$ such that its perimeter is 2: $a+b+c=2$. This equation determines a plane in $abc$-space and, since $a,b,c\geq 0$, the intersection of that plane with the positive orthant contains the parameter space of triangles. This intersection is a simplex, but not every point in this simplex determines a triangle: the side lengths of a triangle must also satisfy the three triangle inequalities
\[
	a \leq b+c, \quad b \leq c + a, \quad c \leq a + b.
\]

\begin{figure}[htbp]
	\centering
		\includegraphics[height=2in]{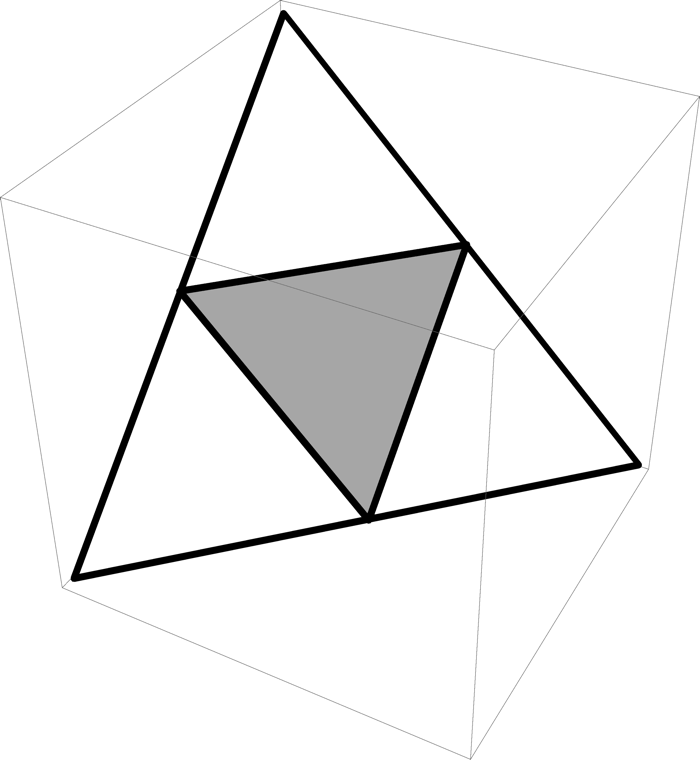}
	\caption{The simplex $a+b+c=2$ in $abc$-space. The darker sub-triangle consists of those points for which $a,b,c$ satisfy the triangle inequalities.}
	\label{fig:simplex}
\end{figure}

Rather than incorporate these slightly ungainly inequalities, we change coordinates, defining
\begin{equation}\label{eqn:sa coords}
	s_a = \frac{b+c-a}{2}, \quad s_b = \frac{c+a-b}{2}, \quad s_c = \frac{a+b-c}{2}.
\end{equation}
These quantities have a long history in triangle geometry, most notably as the radii of three mutually tangent circles centered at the vertices of the triangle.

Since we've already fixed $a+b+c=2$, these coordinates can also be rewritten as $s_a=1-a$, $s_b =1-b$, $s_c=1-c$.\footnote{Here the 1 should be thought of as the semiperimeter of the triangle. This is why we fixed the perimeter to be 2.} Now we see that $s_a+s_b+s_c=3-(a+b+c)=1$ and the triangle inequalities become the conditions $s_a\geq 0$, $s_b \geq 0$, $s_c \geq 0$, so triangle space is parametrized by the standard simplex in $s_a s_b s_c$-space. If we were to perform the analysis from Section~\ref{sec:obtuse} on this simplex, we would get answers similar to those given in the last section of Robin's paper~\cite{Robin:2009vx}.\footnote{Robin is working in the projection of the simplex to the plane, rather than in the simplex itself. This projection distorts the Riemannian metric, so computing in the simplex will produce a slightly different answer than Robin's.}

However, the simplex only has a finite symmetry group, whereas we would prefer, following Portnoy~\cite{Portnoy:1994vh}, a transitive symmetry group, since this induces a canonical probability measure on the parameter space. Therefore, it is desirable to symmetrize by taking square roots: define new coordinates $x,y,z$ by 
\begin{equation}\label{eqn:xyz coords}
	x^2 = s_a, \quad y^2 = s_b, \quad z^2 = s_c,
\end{equation}
so that $x^2+y^2+z^2=s_a+s_b+s_c=1$, and hence the points $(x,y,z)$ lie on the unit sphere. The sphere has a transitive group of symmetries, namely the rotation group, and we actually get rather more than just a canonical probability measure: the round metric on the sphere is an invariant Riemannian metric which is unique up to scale.

Since each of the eight points $(\pm x, \pm y, \pm z)$ map to the same $(s_a,s_b,s_c)$, and hence to the same triangle, the unit sphere is (generically) an eightfold covering of triangle space. The intersection of the sphere with each closed orthant contains points mapping to all possible triangles up to similarity.

Combining \eqref{eqn:sa coords} and \eqref{eqn:xyz coords}, we can translate directly between side lengths $a,b,c$ and sphere coordinates $x,y,z$ using
\begin{equation}\label{eqn:coord translation}
	a = 1-x^2, \quad b = 1-y^2, \quad c = 1-z^2.
\end{equation}

\begin{figure}[htbp]
	\centering
		\includegraphics[height=2in]{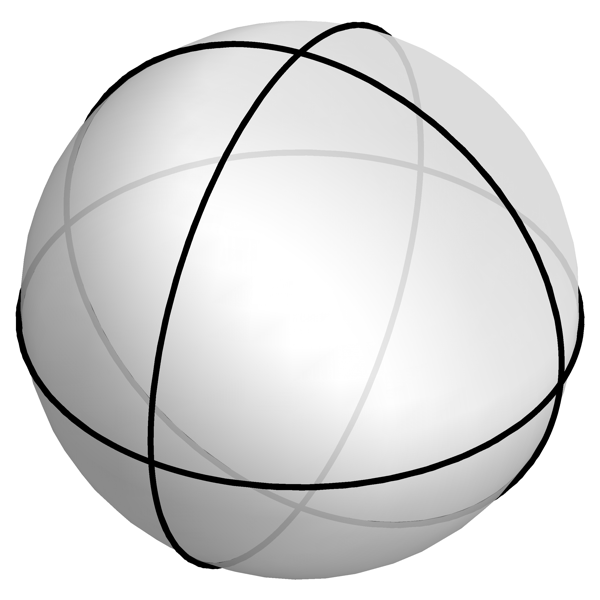}
		\hspace{1in}
		\includegraphics[height=2in]{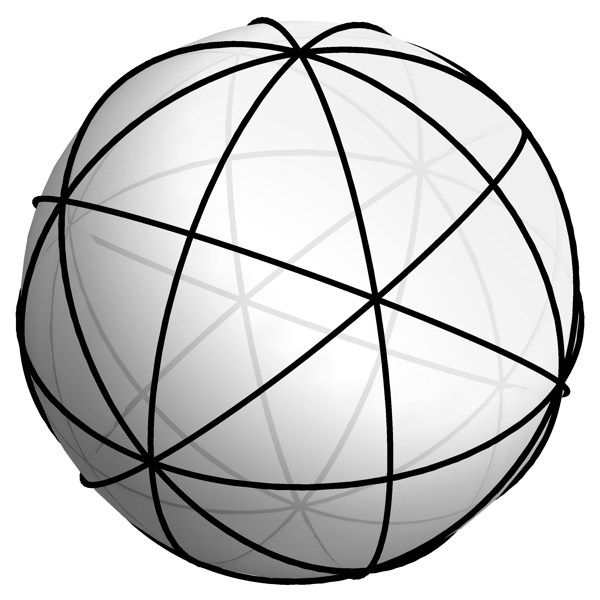}
	\caption{On the left we see the partition of the sphere into the 8 fundamental domains of the group $(\Z/2\Z)^3$ which independently flips the signs of the coordinates. On the right is the partition of the sphere into the 48 fundamental domains of the hyperoctahedral group $B_3 \simeq (\Z/2\Z)^3 \rtimes S_3$ which permutes the coordinates and changes their signs.}
	\label{fig:fundamental domains}
\end{figure}

Using $(a,b,c)$ coordinates, we are implicitly parametrizing \emph{ordered} triangles, where the order of the side lengths matters. If we prefer to think of unordered triangles, we can divide by the action of the permutation group $S_3$, which acts by permuting $a$, $b$, and $c$. This action lifts to the standard permutation action on the sphere, and indeed fits nicely together with the action of changing signs of the coordinates above: the \emph{hyperoctahedral group} $B_3$ acts on the sphere by signed permutations, permuting the coordinates $(x,y,z)$ and changing their signs. $B_3$ is the semidirect product $(\Z/2\Z)^3 \rtimes S_3$ of the group $(\Z/2\Z)^3$ of order $2^3=8$ which acts by changing signs and the permutation group $S_3$ of order $3!=6$. We can also see $B_3$ as the finite subgroup of the isometry group $O(3)$ consisting of all orthogonal matrices with integer entries.

The action of the hyperoctahedral group on polygon space is described in some detail in~\cite{random-triangles}; here we confine ourselves to the observation that $B_3$ acts freely on points of the sphere with all three coordinates having distinct, nonzero magnitudes and it has order 
\[
	|B_3| = 2^3 \times 3! = 48.
\]
Hence, since elements of $B_3$ are isometries, the action of $B_3$ naturally divides the sphere into 48 congruent chambers whose boundaries are the great circles where either two coordinates agree (up to sign) or one coordinate is zero. As seen in Figure~\ref{fig:fundamental domains}, each chamber is a 45--60--90 spherical triangle and we will soon see that we can interpret a chamber as a parameter space for the collection of unordered triangles.

\section{A first solution} \label{sec:first solution}

Given our identification of triangles with the unit sphere, the problem of finding the least symmetric triangle is now simple, at least conceptually: the isosceles triangles are exactly those triangles fixed by some permutation of the edge lengths, so we should identify the subset of the sphere corresponding to isosceles triangles, and then determine the point(s) which are furthest from this subset.

In terms of side length coordinates $a,b,c$, a triangle is isosceles if and only if it satisfies one of the equations
\[
	a = b, \quad b = c, \quad c = a.
\]
Using \eqref{eqn:coord translation} and simplifying slightly, this translates to the equations
\[
	x^2 = y^2, \quad y^2 = z^2, \quad z^2 = x^2
\]
on the sphere. In other words, the subset of the sphere corresponding to isosceles triangles is the intersection of the sphere with the planes
\[
	x \pm y = 0, \quad y \pm z = 0, \quad z \pm x = 0.
\]
Since these are planes through the origin, their intersections with the sphere are great circles, which are geodesics on the sphere. Note that these great circles are a subset of the circles giving the boundaries of the chambers induced by the action of the hyperoctahedral group.

\begin{figure}[htbp]
	\centering
		\includegraphics[height=2in]{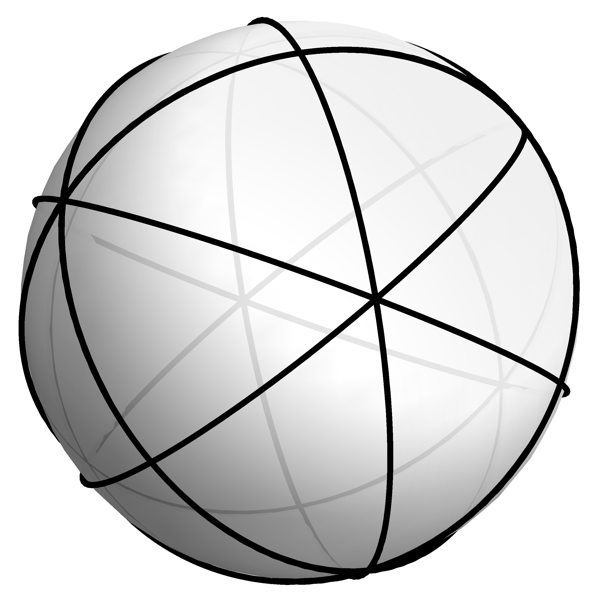}
	\caption{Tiling of the sphere determined by the isosceles triangles.}
	\label{fig:isosceles}
\end{figure}

As seen in Figure~\ref{fig:isosceles}, these great circles determine a tiling of the sphere into 24 spherical triangles, each of which is a 60--60--90 triangle. Therefore, there will be exactly 24 points on the sphere which are furthest from the subset of isosceles triangles: each is the incenter of one of the 24 triangles in the tiling. The corresponding most scalene triangles will all be equivalent up to relabeling the edges, so we can choose any of the 24 triangles and find its incenter.

For example, the curves $x-y=0$, $y+z=0$, and $y-z=0$ determine two of the triangles in the tiling: one with $x$ and $y$ nonnegative the other its antipodal image. We will focus on the triangle $\mathcal{D}$ with $x$ and $y$ both nonnegative, which can also be described as those points $(x,y,z)$ on the sphere satisfying $|z| \leq y \leq x$. Translating to side lengths, these inequalities become $a \leq b \leq c$, so $\mathcal{D}$ parametrizes triangles with side length written in ascending order.

Just as in the plane, the incenter of a spherical triangle is the intersection of the three angle bisectors. The great circles $y+z=0$ and $y-z=0$ are perpendicular and the angle bisector lies on the equator $z=0$, so the incenter will be a point of the form $(x,y,0)$ with $0\leq y \leq x$. 

\begin{figure}[htbp]
	\centering
		\includegraphics[height=2in]{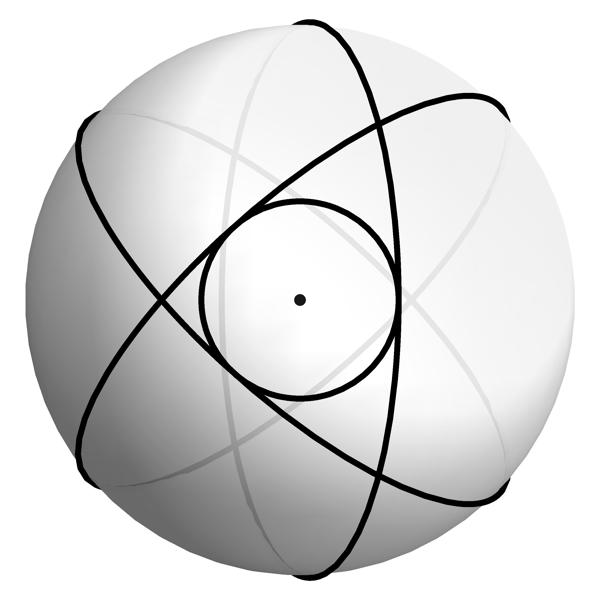}
	\caption{The point $\left(\frac{2}{\sqrt{5}},\frac{1}{\sqrt{5}},0\right)$ corresponding to the least symmetric triangle. It is the incenter of the displayed spherical triangle $\mathcal{D}$. The radius (along the sphere) of the inscribed circle is $\arccos\left(\frac{3}{\sqrt{10}}\right) \approx 0.32175$.}
	\label{fig:inscribed circle}
\end{figure}

In order to determine the angle bisector of $x-y=0$ and $y-z=0$, we will consider the unit normal vectors to the planes, namely $v_1 = \left(-\frac{1}{\sqrt{2}}, \frac{1}{\sqrt{2}},0\right)$ and $v_2 = \left(0,\frac{1}{\sqrt{2}}, -\frac{1}{\sqrt{2}}\right)$. These vectors form an angle
\[
	\theta = \arccos(v_1 \cdot v_2) = \arccos(1/2) = \pi/3
\]
and the vector halfway between them is proportional to the sum
\[
	v_1 + v_2 = \left(-\frac{1}{\sqrt{2}}, \frac{1}{\sqrt{2}},0\right)+ \left(0,\frac{1}{\sqrt{2}}, -\frac{1}{\sqrt{2}}\right) = \left(-\frac{1}{\sqrt{2}},\sqrt{2},-\frac{1}{\sqrt{2}}\right).
\]
Scaling the above vector by $\sqrt{2}$, we see that the great circle which bisects $x-y=0$ and $y-z=0$ is
\[
	-x+2y-z=0.
\]
The intersection of this great circle with the equator $z=0$ is our desired point $\left(\frac{2}{\sqrt{5}},\frac{1}{\sqrt{5}},0\right)$ (see Figure~\ref{fig:inscribed circle}), corresponding to the triangle with side lengths
\begin{align*}
	a & = 1-\left(\frac{2}{\sqrt{5}}\right)^2 = \frac{1}{5} \\
	b & = 1-\left(\frac{1}{\sqrt{5}}\right)^2 = \frac{4}{5} \\
	c & = 1-0^2 = 1.
\end{align*}

As a proposed least symmetric triangle, the fact that the side lengths are in a $1 : 4 : 5$ ratio is gratifying; rather less so is that this is a degenerate triangle with all three sides lying in a line. It is, perhaps, not surprising that the least symmetric triangle would have as much difference as possible between its shortest and longest side lengths subject to the constraints of the triangle inequalities, resulting in a length ratio of $1 : r : r+1$. Both as $r \to 1$ and as $r \to \infty$, the resulting triangle becomes isosceles, so the specific value $r=4$ is apparently a balance between these two extremes subject to the constraints of the geometry of the sphere.

However, degenerate triangles like the one we just found are actually symmetric: reflecting across the line of degeneracy fixes the degenerate triangle. The issue is that $\mathcal{D}$ is, in fact, (almost) a double covering of the space of unordered triangles. The points $(x,y,z)$ and $(x,y,-z)$ map to the same triangle, so almost all triangles have two preimages in $\mathcal{D}$; the exceptions are those triangles with $1 = c = 1-z^2$, which come from points of the form $(x,y,0)$. Therefore, the set of triangles we consider symmetric should include not just the isosceles triangles, but also the degnerate triangles.

\section{Excluding degenerate triangles}\label{sec:non-degenerate}

The degenerate triangles are those with one side length being half the total perimeter of the triangle: this forces the other two sides to lie in the same line as the long side. Given our normalization that triangles should have perimeter 2, this means that the degenerate triangles are those with $a=1$, $b=1$, or $c=1$.\footnote{Note that a triangle can also be doubly-degenerate, with two sides of length 1 and one of length 0. This means that two vertices coincide.}

In terms of $x,y,z$ coordinates, the degenerate triangles are those with
\[
	x=0 \quad \text{or} \quad \quad y=0 \quad \text{or} \quad z=0.
\]
Adding these three great circles to the six corresponding to the isosceles triangles gives the tiling of the sphere by 48 congruent 45--60--90 triangles that we saw in Figure~\ref{fig:fundamental domains}, and we have now justified the claim that each triangle can be thought of as the space of unordered triangles up to similarity. 

\begin{figure}[htbp]
	\centering
		\includegraphics[height=2in]{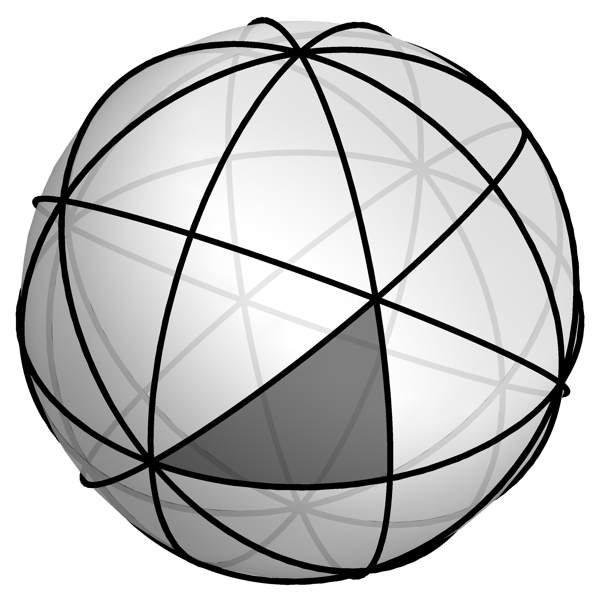}
	\caption{The triangle $\mathcal{T}$ bounded by $x-y=0$, $y-z=0$, and $z=0$, the interior of which parametrizes scalene triangles.}
	\label{fig:tiling2}
\end{figure}

For specificity, we will focus on the spherical triangle $\mathcal{T}$ bounded by $x-y=0$, $y-z=0$, and $z=0$ as our model of the space of triangles. Equivalently, $\mathcal{T}$ is the subset of the sphere with $0 \leq z \leq y \leq x$, which means that the side lengths $a,b,c$ of the triangles corresponding to points in $\mathcal{T}$ satisfy
\[
	0 \leq a \leq b \leq c \leq 1,
\]
and the interior of $\mathcal{T}$ corresponds to those triangles for which all of the above inequalities are strict.

\begin{figure}[htbp]
	\centering
		\includegraphics[height=2in]{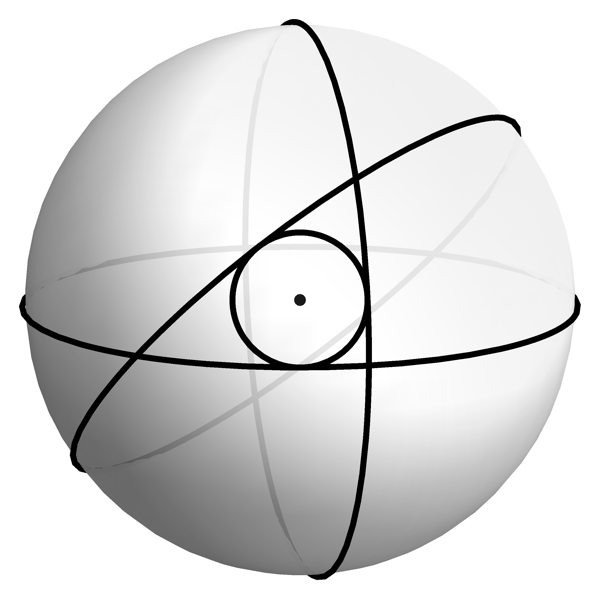}
	\caption{The incenter and incircle of the spherical triangle $\mathcal{T}$ which parametrizes the non-degenerate scalene triangles. The incenter is the point $\frac{1}{\sqrt{13+6 \sqrt{2}}} \left(1+2 \sqrt{2},1+\sqrt{2},1\right)$ and the radius of the inscribed circle is $\arcsin\left(\sqrt{\frac{1}{97} \left(13-6 \sqrt{2}\right)}\right) \approx 0.217449$.}
	\label{fig:incenter2}
\end{figure}

Therefore, a more plausible ``least symmetric triangle'' will be the incenter of this triangle, as visualized in Figure~\ref{fig:incenter2}. Since this will again be the intersection of the angle bisectors, it must be the intersection of the previously determined great circle $-x+2y-z=0$ (which bisects the angle formed by the sides $x-y=0$ and $y-z=0$) and the great circle halfway between $x-y=0$ and $z=0$. 

\begin{figure}[htbp]
	\centering
		\includegraphics[height=1in]{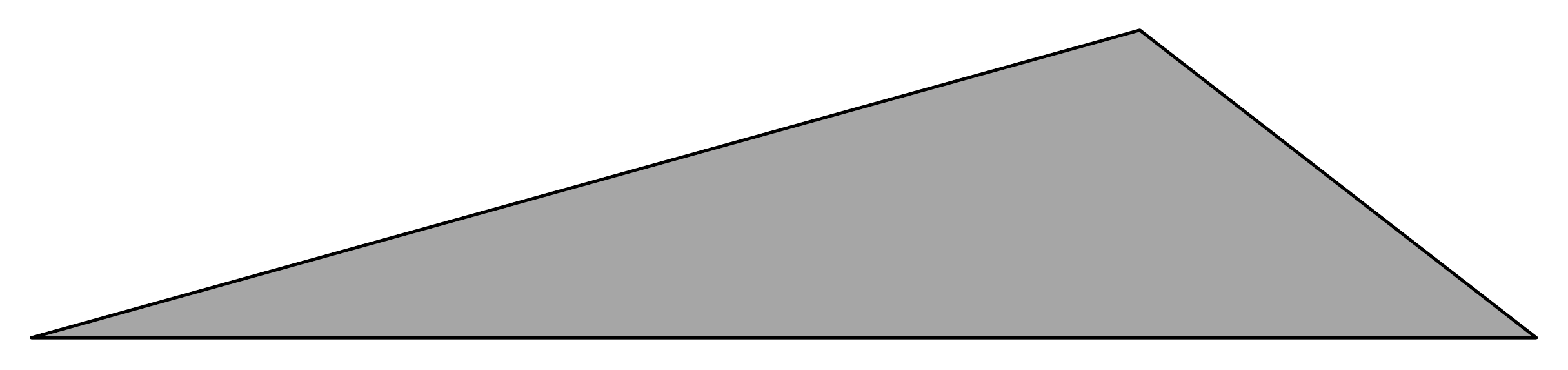}
	\caption{The least symmetric triangle. Its side lengths fall in the ratio $1: 3-\frac{1}{\sqrt{2}} : 3$.}
	\label{fig:most scalene v1}
\end{figure}

As before, we can find the latter great circle as being perpendicular to the sum of the unit normal vectors
\[
	\left(-\frac{1}{\sqrt{2}}, \frac{1}{\sqrt{2}},0\right) + (0,0,1) = \left(-\frac{1}{\sqrt{2}}, \frac{1}{\sqrt{2}},1\right).
\]
Multiplying by $\sqrt{2}$, we see that we're looking at the great circle $-x+y+\sqrt{2}z=0$, which intersects $-x+2y-z=0$ at the point
\[
	(x,y,z) = \frac{1}{\sqrt{13+6 \sqrt{2}}} \left(1+2 \sqrt{2},1+\sqrt{2},1\right) \approx (0.825943, 0.520841, 0.215739).
\]
We have now proved:

\begin{proposition}\label{prop:nondegenerate}
	The least symmetric triangle has side lengths
\begin{align*}
	a & = 1-x^2 = \frac{28+2\sqrt{2}}{97} \approx 0.3178 \\
	b & = 1-y^2 = \frac{82-8\sqrt{2}}{97} \approx 0.7287 \\
	c & = 1-z^2 = \frac{84+6\sqrt{2}}{97} \approx 0.9524
\end{align*}
and is shown in Figure~\ref{fig:most scalene v1}. 
\end{proposition}

Though not quite as simple as for the degenerate triangle from Section~\ref{sec:first solution}, the side lengths of this triangle still form the pleasing ratio $1: 3-\frac{1}{\sqrt{2}} : 3$.

\section{Refinements}
\label{sec:obtuse}

Having found the least symmetric triangle, we can hardly resist refining the question and trying to find the least symmetric obtuse and acute triangles, even though \emph{obtuse} and \emph{acute} are not notions that naturally generalize to $n$-gons. To do so, we want to find the point on the sphere furthest from not only the isosceles and the degenerate triangles, but also from the right triangles. Hence, we need to identify the collection of points on the sphere corresponding to the right triangles. 

In terms of the side lengths $a,b,c$, the right triangles are uniquely characterized by satisfying the Pythagorean theorem, meaning that
\[
	a^2 + b^2 = c^2 \quad \text{or} \quad b^2 + c^2 = a^2 \quad \text{or} \quad c^2 + a^2 = b^2.
\]
Translating into $x,y,z$ coordinates, the subset of the sphere corresponding to right triangles is the set of points satisfying one of the following quartic equations:
\[
	(1-x^2)^2 + (1-y^2)^2 = (1-z^2)^2, \  (1-y^2)^2 + (1-z^2)^2 = (1-x^2)^2, \  (1-z^2)^2 + (1-x^2)^2 = (1-y^2)^2.
\]

\begin{figure}[htbp]
	\centering
		\includegraphics[height=2in]{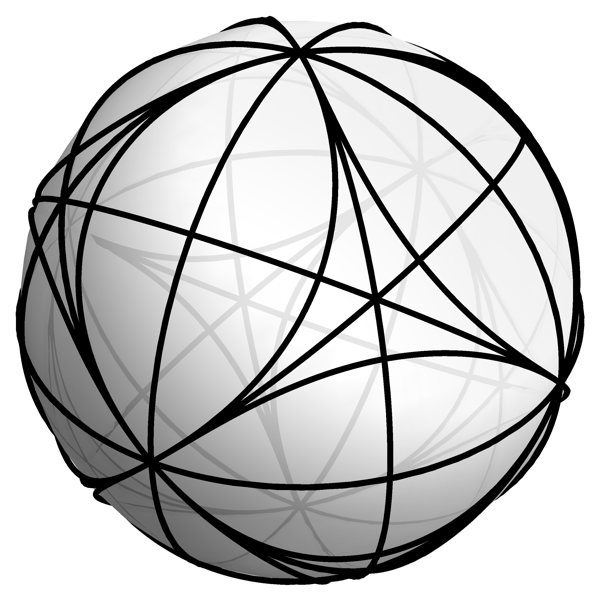}
		\hspace{1in}
		\includegraphics[height=2in]{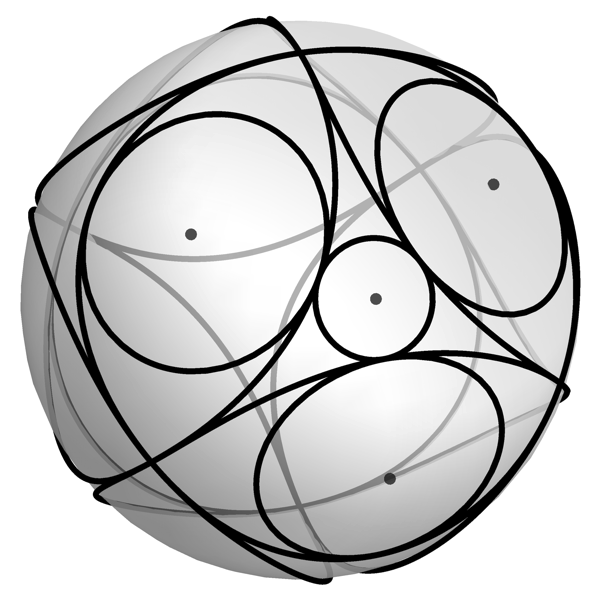}
	\caption{The curves on the left represent the degenerate triangles, the isosceles triangles, and the right triangles. On the right we see the points corresponding to the equilateral triangle and to the degenerate isosceles triangles, which (locally) maximize distance from the right triangles. In both cases the closest right triangle is a 45--45--90 triangle, which is at a distance $\arccos\left(\frac{\sqrt{2}-1+2 \sqrt{\sqrt{2}-1}}{\sqrt{3}}\right)\approx 0.188401$ from the equilateral triangle and $\arccos \left(\sqrt{2 \left(\sqrt{2}-1\right)}\right)\approx 0.427079$ from the degenerate isosceles triangles.}
	\label{fig:tiling3}
\end{figure}

Incidentally, as in Figure~\ref{fig:tiling3}, the most acute triangle (in the sense of being furthest from the right triangles) is the equilateral triangle, and the most obtuse is the degenerate triangle with side lengths $(1/2,1/2,1)$. This is reassuringly intuitive: the equilateral triangle is surely the acute triangle least like a right triangle, and we would expect the most obtuse triangle to have a $180^\circ$ angle.

\subsection{Obtuse triangles} 
\label{sub:obtuse}

To find the least symmetric obtuse triangle, we again focus our attention on the spherical triangle $\mathcal{T}$ bounded by $x-y=0$, $y-z=0$, and $z=0$. Since the curve of right triangles satisfying the equation $(1-x^2)^2+(1-y^2)^2 = (1-z^2)^2$ is the one which intersects $\mathcal{T}$, we focus on it. Since this curve is not a geodesic, the region $\mathcal{O}$ bounded by it and the great circles $x-y=0$ and $z=0$ is not a spherical triangle, making it a more substantial challenge to find the point which is maximally far from the boundary of $\mathcal{O}$. 

The point we are after must be equidistant from the great circles $x-y=0$ and $z=0$, so it must lie on the great circle which bisects the angle between them, namely the great circle $-x+y+\sqrt{2}z=0$ that we already found in Section~\ref{sec:non-degenerate}. This great circle contains the orthonormal vectors 
\[
	u_1 = \left(\frac{1}{\sqrt{2}},\frac{1}{\sqrt{2}}, 0\right), \quad u_2 = \left(\frac{1}{2},-\frac{1}{2},\frac{1}{\sqrt{2}}\right),
\]
so it can be parametrized as
\[
	p(t) = \cos t\, u_1 + \sin t\, u_2 = \left(\frac{\sqrt{2}\cos t + \sin t}{2}, \frac{\sqrt{2} \cos t - \sin t}{2}, \frac{\sin t}{\sqrt{2}}\right).
\]

The angle that $p(t)$ makes with the $z$-axis is simply 
\[
	\theta(t) = \arccos(p(t) \cdot (0,0,1)) = \arccos\left(\frac{\sin t}{\sqrt{2}}\right),
\]
which means that the spherical distance from $p(t)$ to the great circle $z=0$ (and hence also to the great circle $x-y=0$) is
\[
	\frac{\pi}{2} - \theta(t) = \frac{\pi}{2} - \arccos\left(\frac{\sin t}{\sqrt{2}}\right) = \arcsin\left(\frac{\sin t}{\sqrt{2}}\right)
\] 
since $\sin \phi = \cos(\pi/2-\phi)$ for any $\phi$.

Therefore, the problem is to determine the value of $t$ for which the spherical distance from $p(t)$ to the curve $(1-x^2)^2+(1-y^2)^2 = (1-z^2)^2$ is equal to $ \arcsin\left(\frac{\sin t}{\sqrt{2}}\right)$.

Using the fact that our points are on the sphere $x^2+y^2+z^2=1$, we can parametrize the curve of right triangles $(1-x^2)^2+(1-y^2)^2 = (1-z^2)^2$ by
\[
	q(x) = \left(x, \sqrt{\frac{1-x^2}{1+x^2}},x \sqrt{\frac{1-x^2}{1+x^2}}\right),
\]
so the spherical distance from $p(t)$ to $q(x)$ is
\begin{multline*}
	d(t,x)=\arccos(p(t)\cdot q(x)) = \arccos \left(x \frac{\sqrt{2}\cos t + \sin t}{2} + \sqrt{\frac{1-x^2}{1+x^2}}\frac{\sqrt{2} \cos t - \sin t}{2}\right.\\
	\left.+x \sqrt{\frac{1-x^2}{1+x^2}}\frac{\sin t}{\sqrt{2}}\right).
\end{multline*}

Therefore, we are seeking the smallest positive value of $t$ so that 
\[
	d(t,x) = \arcsin\left(\frac{\sin t}{\sqrt{2}}\right)
\]
or, after eliminating the inverse trig functions, so that
\[
	x \frac{\sqrt{2}\cos t + \sin t}{2} + \sqrt{\frac{1-x^2}{1+x^2}}\frac{\sqrt{2} \cos t - \sin t}{2}+x \sqrt{\frac{1-x^2}{1+x^2}}\frac{\sin t}{\sqrt{2}} = \sqrt{1-\frac{\sin^2 t}{2}}.
\]
This equation can be solved for $x$ in, e.g., \emph{Mathematica}, though the solution is not pleasant: printing out all four solutions would require hundreds of pages. In turn, given $x$ as a function of $t$, the challenge is to find the smallest $t$ for which $x$ is real. Since the process of finding it was \emph{ad hoc} and unenlightening, we simply report that the smallest such $t$ is
\[
	t_0 = 2 \arctan \alpha,
\]
where $\alpha \approx 0.140112$ is the smallest positive root of the even, palindromic polynomial
\begin{multline*}
	16 z^{24}-992 z^{22}+9689 z^{20}-36232 z^{18}+100908 z^{16}-197080 z^{14}+238166z^{12}\\
	   -197080 z^{10}+100908 z^8-36232 z^6+9689 z^4-992 z^2+16.
\end{multline*}
Complete details can be found in the supplementary materials~\cite{supplementary-materials}.

\begin{figure}[htbp]
	\centering
		\includegraphics[height=2in]{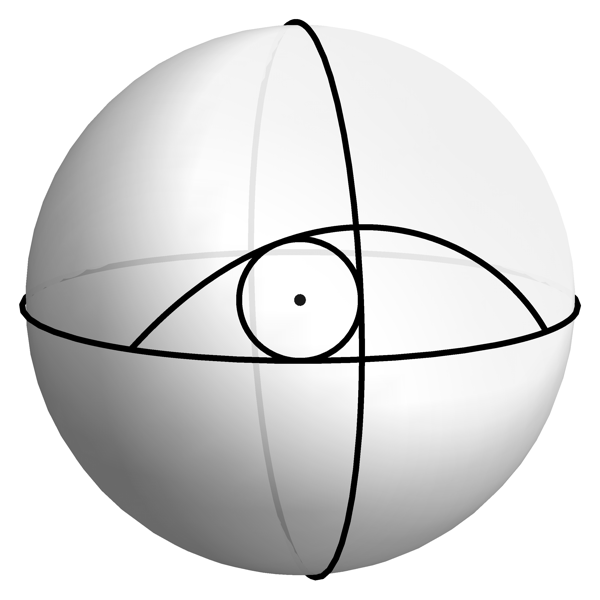}
		\hspace{1in}
		\includegraphics[height=2in]{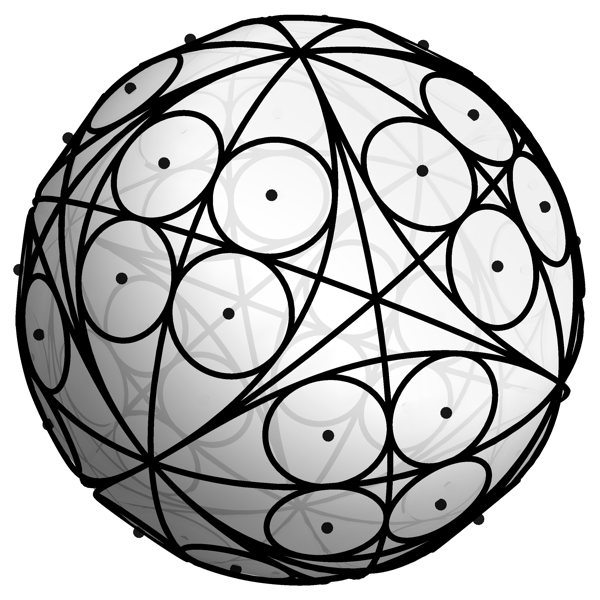}
	\caption{The figure on the left shows the region $\mathcal{O}$ together with the point $p(t_0)$ maximally far from the boundary of $\mathcal{O}$ and the circle of radius $\arcsin\left(\frac{\sqrt{2} \alpha }{1+\alpha ^2}\right)\approx 0.195578$ around the point. The figure on the right shows the 48 different points and circles under the action of the hyperoctahedral group $B_3$ which permutes the coordinates and independently changes their signs.}
	\label{fig:inscribed3}
\end{figure}

\begin{figure}[htbp]
	\centering
		\includegraphics[height=1.3in]{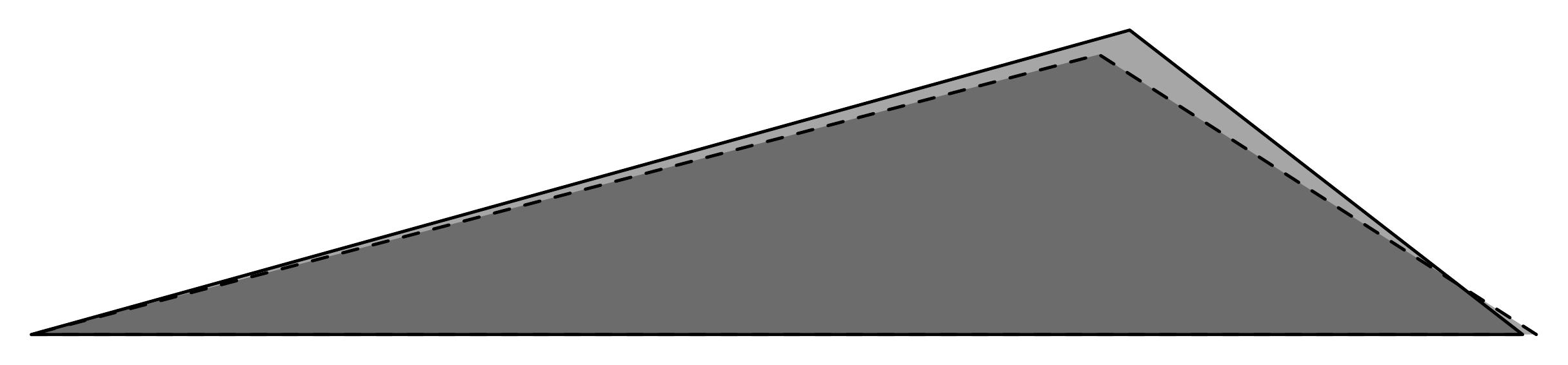}
	\caption{A comparison between asymmetric triangles. Not surprisingly, the (dashed) least symmetric obtuse triangle has a larger obtuse angle than the (solid) least symmetric triangle.}
	\label{fig:mostscalene2}
\end{figure}

Therefore, the point on the sphere that we're after is 
\[
	p(t_0) \approx (0.817293, 0.542464, 0.194334),
\]
and converting to side length coordinates yields:

\begin{proposition}\label{prop:obtuse}
	The least symmetric obtuse triangle has side lengths
	\begin{multline*}
		\frac{1}{2(1+\alpha^2)^2}\left(1-2\sqrt{2}\alpha+4\alpha^2+2\sqrt{2}\alpha^3+\alpha^4,1+2\sqrt{2}\alpha+4\alpha^2-2\sqrt{2}\alpha^3+\alpha^4,2+2\alpha^4\right) \\
		\approx (0.332032,0.705733,0.962234).
	\end{multline*}
\end{proposition}


\subsection{Acute triangles} 
\label{sub:acute}
\begin{figure}[htbp]
	\centering
		\includegraphics[height=2in]{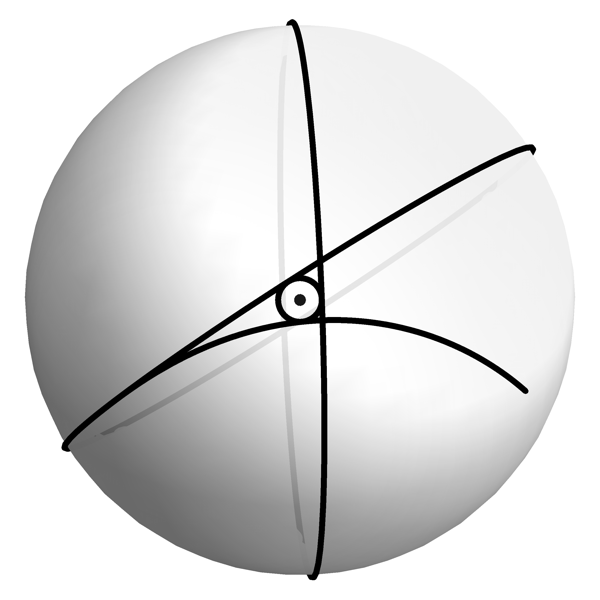}
		\hspace{1in}
		\includegraphics[height=2in]{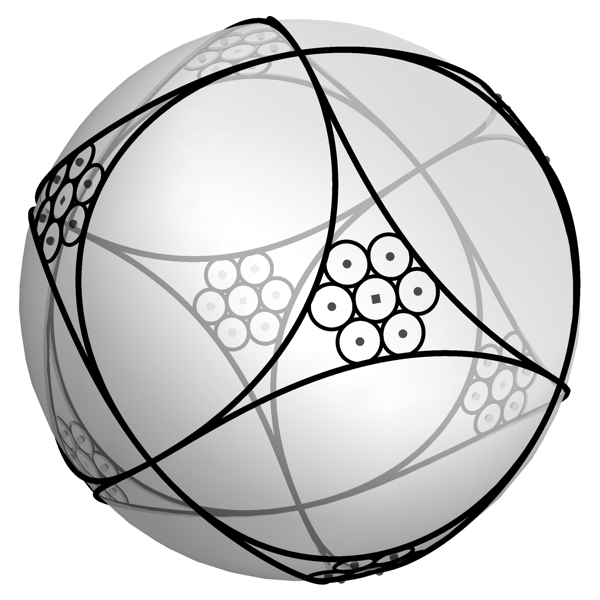}
	\caption{The point $\widetilde{p}(\widetilde{t}_0) \approx (0.670125, 0.571734, 0.473343)$ in $\mathcal{A}$ maximally far from the curves corresponding to the isosceles and right triangles. The common distance to the boundary of the region $\mathcal{A}$ is $\arcsin \left(\frac{\widetilde{\alpha}}{1+\widetilde{\alpha}^2}\right)\approx 0.069629$. Permuting the coordinates of $\widetilde{p}(\widetilde{t}_0)$ yields a ring of six points equidistant from the equilateral triangle point $\left(\frac{1}{\sqrt{3}},\frac{1}{\sqrt{3}},\frac{1}{\sqrt{3}}\right)$. The difference is too small to see, but the radius of the circle around the equilateral point ({\tiny $\blacksquare$}) in the right figure is smaller than that of the inscribed circle around $\widetilde{p}(\widetilde{t}_0)$ by $\approx 0.000333$.}
	\label{fig:inscribed5}
\end{figure}

Now, we turn to solving the corresponding problem for acute triangles. The subset $\mathcal{A} \subseteq \mathcal{T}$ corresponding to acute triangles is shown on the left in Figure~\ref{fig:inscribed5}; it is bounded by the great circles $x-y=0$ and $y-z=0$, as well as the curve of right triangles $(1-x^2)^2+(1-y^2)^2 = (1-z^2)^2$. The point which maximizes distance to the boundary must lie on the great circle halfway between $x-y=0$ and $y-z=0$, which has equation $x-2y+z=0$ and can be parametrized by
\[
	\widetilde{p}(t) = \cos t \left(\frac{1}{\sqrt{3}},\frac{1}{\sqrt{3}},\frac{1}{\sqrt{3}}\right) + \sin t \left(\frac{1}{\sqrt{2}},0,-\frac{1}{\sqrt{2}}\right).
\]

The distance from $\widetilde{p}(t)$ to the boundary great circle $x-y=0$ is $\arcsin\left(\frac{\sin t}{2}\right)$, and the distance to the point $q(x)$ on the curve of right triangles is
\[
	\widetilde{d}(t,x) = \arccos\left(\frac{\sqrt{\frac{1-x^2}{x^2+1}} \cos t}{\sqrt{3}}+x \sqrt{\frac{1-x^2}{x^2+1}} \left(\frac{\cos t}{\sqrt{3}}-\frac{\sin t}{\sqrt{2}}\right)+x \left(\frac{\sin t}{\sqrt{2}}+\frac{\cos t}{\sqrt{3}}\right)\right).
\]
Hence, the challenge is to determine the smallest value of $t$ for which $\widetilde{d}(t,x) = \arcsin\left(\frac{\sin t}{2}\right)$. Again, we solve this equation for $x$, and then find the smallest $t$ which makes $x$ real, which is
\[
	\widetilde{t}_0 = 2 \arctan \widetilde{\alpha}
\]
where $\widetilde{\alpha} \approx 0.069912$ is the smallest positive root of the even, palindromic polynomial
\begin{align*}
	& 131072 z^{52}-30081024 z^{50}+715784192 z^{48}-10181738496 z^{46}+83609604096 z^{44}\\
	& -443259328512 z^{42}+1410471953408 z^{40}-1858643071488 z^{38}+18137673285920 z^{36}\\
	& -14367112128688 z^{34}+56162265469488 z^{32}-73041229883512 z^{30}+73382345772378 z^{28}\\
	& -122601623733111 z^{26}+73382345772378 z^{24}-73041229883512 z^{22}+56162265469488 z^{20}\\
	& -14367112128688 z^{18}+18137673285920 z^{16}-1858643071488 z^{14}+1410471953408 z^{12}\\
	& -443259328512 z^{10}+83609604096 z^8-10181738496 z^6+715784192 z^4-30081024 z^2+131072.
\end{align*}
Again, more details are in the supplementary materials~\cite{supplementary-materials}. 

Therefore, the incenter of $\mathcal{A}$ is
\[
	\widetilde{p}(\widetilde{t}_0) \approx (0.670125, 0.571734, 0.473343)
\]
and we can convert to side lengths to conclude:

\begin{proposition}\label{prop:acute}
	The least symmetric acute triangle has side lengths
	\begin{multline*}
		\frac{2}{3(1+\widetilde{\alpha}^2)^2}\left(1-\sqrt{6}\widetilde{\alpha}+\widetilde{\alpha}^2+\sqrt{6}\widetilde{\alpha}^3+\widetilde{\alpha}^4,1+4\widetilde{\alpha}^2+\widetilde{\alpha}^4,1+\sqrt{6}\widetilde{\alpha}+\widetilde{\alpha}^2-\sqrt{6}\widetilde{\alpha}^3+\widetilde{\alpha}^4\right)\\
		\approx (0.550933, 0.673120, 0.775946).
	\end{multline*}
\end{proposition}

\begin{figure}[htbp]
	\centering
		\includegraphics[height=1.4in]{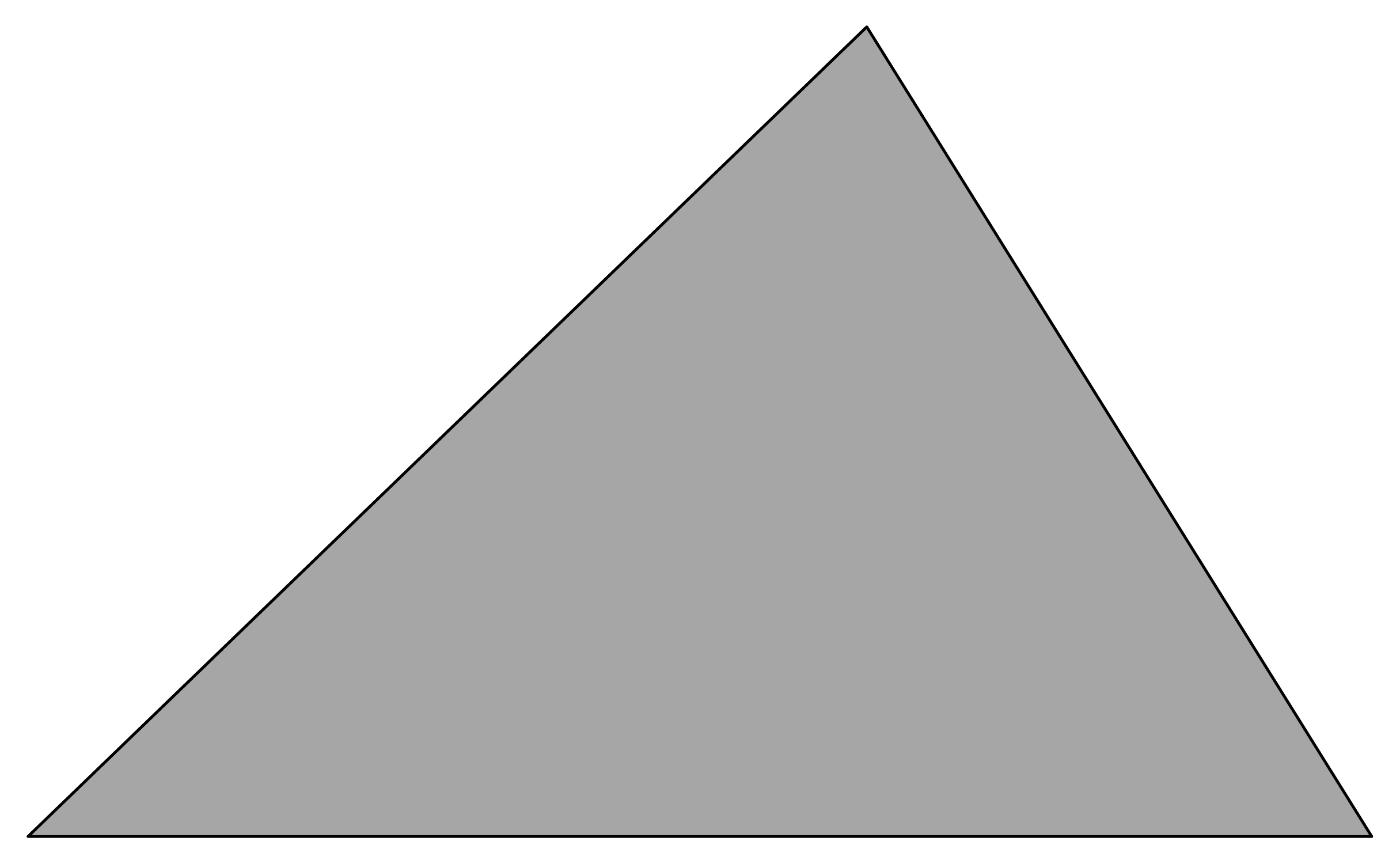}
	\caption{The least symmetric acute triangle, with side lengths $\approx (0.550933, 0.67312, 0.775946)$.}
	\label{fig:mostacute}
\end{figure}


\section{Conclusion and open questions}

The emphasis in the chemistry literature is on finding \emph{chirality measures} whose maxima are considered the most chiral triangles. Of course, our least symmetric triangle from Proposition~\ref{prop:nondegenerate} can also be thought of in this way: it maximizes the minimum distance to the 9 great circles representing isosceles and degenerate triangles. Since this function is continuous but only piecewise smooth, a natural question to ask is: does there exist a reasonable smooth function on the sphere which vanishes precisely at the isosceles and degenerate triangles and which is maximized on the hyperoctahedral group orbit of the point $\frac{1}{\sqrt{13+6 \sqrt{2}}} \left(1+2 \sqrt{2},1+\sqrt{2},1\right)$?\footnote{Rassat and Fowler~\cite{Rassat:2003jn} showed that \emph{any} non-isosceles triangle is the most chiral triangle according to some chirality measure in a particular family, but their chirality measures are slightly unnatural to define on the sphere.} And similarly for the other special points we have found? At the very least, we expect, based on numerical experiments and O'Hara's results on planar convex bodies~\cite{OHara:2012fm}, that these points are the limits of maxima of suitably renormalized Riesz-type potentials as the exponent goes to $-\infty$.

The problem of extending this analysis to $n$-gons is substantial and interesting. As alluded to in the introduction, the spherical tiling induced by the hyperoctahedral group action on the sphere (Figure~\ref{fig:fundamental domains}) generalizes to a decomposition of $n$-gon space for any $n$, namely the images of the fundamental domain of the standard hyperoctahedral group action on $G_2(\R^n)$. In Section~\ref{sec:non-degenerate} we found the incenter of $\mathcal{T}$ and hence, by taking this point's hyperoctahedral group orbit, of all the other triangles in the tiling. What are the incenters of the cells in the decomposition of $n$-gon space, which we can interpret as the least symmetric $n$-gon?

Similarly, what is the least symmetric $n$-gon in $\R^3$? Just as $n$-gons in the plane are modeled by $G_2(\R^n)$, $n$-gons in space correspond to points in the complex Grassmannian $G_2(\C^n)$. The natural generalization of the hyperoctahedral group is given by replacing $(\Z/2)^n = (O(1))^n$ with $(U(1))^n$, yielding the group $(U(1))^n \rtimes S_n$, which acts on $G_2(\C^n)$. The fundamental domain (in the sense of Hermann~\cite{Hermann:1962eu}) of this action gives the space of unordered $n$-gons in $\R^3$, and it seems challenging to describe this domain and then to find the point furthest from the boundary.

Finding the least symmetric triangle can be interpreted as finding the optimal shape satisfying certain constraints. As in Section~\ref{sec:obtuse}, where we added the constraints that triangles should be obtuse or acute, other constraints are also interesting; for example, what is the most knotted\footnote{Meaning furthest from the subset of unknotted polygons.} trefoil knot in the space of $23$-gons? Since our polygon model has generalizations to continuous curves in the plane~\cite{Younes:2008gy} and to framed curves in space~\cite{Needham:2017vn}, we are not restricted to asking such questions only about polygons.

Finally, while we have presented the triangle from Proposition~\ref{prop:nondegenerate} as the least symmetric triangle, we could also think of it as the \emph{median} triangle since it is equidistant from the three distinct pieces of the boundary of the region $\mathcal{T}$ parametrizing unordered triangles. This framing suggests the obvious question: what is the \emph{mean} triangle? Or, indeed, the mean $n$-gon? We intend to address this question in future work.

\section*{Acknowledgments}
We would like to thank Vance Blankers, Jason Cantarella, Renzo Cavalieri, Andy Fry, Tom Needham, Eric Rawdon, and Gavin Stewart for stimulating conversations about the geometry of polygon space. We are especially grateful to Noah Otterstetter for his participation in our early conversations about this project.

This work was supported by a grant from the Simons Foundation (\#354225, CS).

\bibliography{scalene-special,papers-export}

\end{document}